\documentclass[aos,MSNbibl,nameyear,seceqn,dvips]{arximspdf}

\doi{10.1214/12-AOS992} 
\volume{40}
\issue{3}
\pubyear{2012}
\firstpage{1665}
\lastpage{1681}

\makeatletter
\newtheorem{Theorem}{Theorem}
\newtheorem{Proposition}{Proposition}
\newtheorem{Lemma}{Lemma}
\makeatother

\begin{document}
\begin{frontmatter}

\title{Identifying locally optimal designs for nonlinear models: A simple extension with profound consequences}
\runtitle{Locally optimal designs}

\begin{aug}
\author[A]{\fnms{Min} \snm{Yang}\corref{}\thanksref{t1}\ead[label=e1]{myang2@uic.edu}}
\and
\author[B]{\fnms{John} \snm{Stufken}\thanksref{t2}\ead[label=e2]{jstufken@uga.edu}}
\runauthor{M. Yang and J. Stufken}
\thankstext{t1}{Supported by NSF Grants DMS-07-07013 and DMS-07-48409.}
\thankstext{t2}{Supported by NSF Grants DMS-07-06917 and DMS-10-07507.}

\affiliation{University of Illinois at Chicago and University of Georgia}

\address[A]{Department of Mathematics, Statistics,\\
\quad and Computer Sciences\\
University of Illinois at Chicago\\
Chicago, Illinois 60607\\
USA\\
\printead{e1}} 
\address[B]{Department of Statistics\\
The University of Georgia\\
Athens, Georgia 30602-7952\\
USA\\
\printead{e2}}
\end{aug}

\received{\smonth{8} \syear{2011}}
\revised{\smonth{2} \syear{2012}}

%
\begin{abstract}
We extend the approach in [\textit{Ann. Statist.} \textbf{38} (2010) 2499--2524]
for identifying locally optimal designs for nonlinear models.
Conceptually the extension is relatively simple, but the consequences
in terms of applications are profound. As we will demonstrate, we can
obtain results for locally optimal designs under many optimality
criteria and for a larger class of models than has been done hitherto.
In many cases the results lead to optimal designs with the minimal
number of support points.
\end{abstract}

%
\begin{keyword}[class=AMS]
\kwd[Primary ]{62K05}
\kwd[; secondary ]{62J12}.
\end{keyword}
\begin{keyword}
\kwd{Locally optimal}
\kwd{Loewner ordering}
\kwd{principal submatrix}
\kwd{support points}.
\end{keyword}

\end{frontmatter}

\section{Introduction}
During the last decades nonlinear models have become a workhorse for
data analysis in many applications. While there is now an extensive
literature on data analysis for such models, research on design
selection has not kept pace, even though there has seen a spike in
activity in recent years. Identifying optimal designs for nonlinear
models is indeed much more difficult than the much better studied
corresponding problem for linear models. For nonlinear models results
can typically only be obtained on a case-by-case basis, meaning that
each combination of model, optimality criterion and objective of the
experiment requires its own proof.

Another challenge is that for a nonlinear model an optimal design
typically depends on the unknown parameters. This leads to the concept
of locally optimal designs, which are optimal for a priori chosen
values of the parameters. The designs may be poor if the choice of
values is far from the true values. Where feasible, a multistage
approach could help with this. A~small initial design is then used to
obtain some information about the parameters,\vadjust{\goodbreak} and this information is
used at the next stage to estimate the true parameter values and to
extend the initial design in a locally optimal way to a larger design.
The design at this second stage could be the final design, or there
could be additional stages at which more design points are selected.
The solution presented in this paper is applicable for a one-shot
approach for finding a locally optimal design as well as for a
multistage approach. The argument that our method can immediately be
applied for the multistage approach is exactly as in \citet{ys09}.

For a broader discussion on the challenges to identify optimal designs
for generalized linear models, many of which apply also for other
nonlinear models, we refer the reader to \citet{k06}.

The work presented here is an extension of \citet{ys09},
\citet{y09} and \citet{d05}. The analytic approach in those
papers unified and extended many of the results on locally optimal
designs that were available through the so-called geometric approach.
The extension in the current paper has major consequences for two
reasons. First, it enables the application of the basic approach in the
three earlier papers to many models for which it could until now not be
used. As a result, this paper opens the door to finding locally optimal
designs for models where no feasible approach was known so far. Second,
for a number of models for which answers could be obtained by earlier
work, the current extension enables the identification of locally
optimal designs with a smaller support. This is important because it
simplifies the search for optimal designs, whether by computational or
analytical methods. Section~\ref{sec4} will illustrate the impact of
our results.

The basic approach in \citet{ys09}, \citet{y09} and \citet{d05},
which is also adopted here, is to identify a subclass
of designs with a simple format, so that for any given design $\xi$,
there exists a~design $\xi^*$ in that subclass with $I_{\xi^*}\geq
I_{\xi}$ under the Loewner ordering. We will refer to this subclass as
a complete class for this problem. Here, $I_{\xi^*}$ and $I_{\xi}$ are
information matrices for a parameter vector $\theta$ under $\xi^*$ and
$\xi$, respectively. Others, such as \citet{p89} have called such
a class \textit{essentially} complete, which is admittedly indeed more
accurate, but also more cumbersome. When searching for a locally
optimal design, for the common information-based optimality criteria,
including $A$-, $D$-, $E$- and $\Phi_p$-criteria, one can thus restrict
consideration to this complete class, both for a one-shot or multistage
approach. Also, as shown in \citet{ys09}, this conclusion
holds for arbitrary functions of the parameters. Ideally, the same
complete class results would apply for all a priori values of the
parameter vector~$\theta$. However, it turns out, as we will see in
Section~\ref{sec4}, that there are instances where complete class results hold
only for certain a priori values of $\theta$.\looseness=-1

\citet{ys09}, \citet{y09} and \citet{d05}
identify small complete classes for certain models. They do so by
showing that for any design $\xi$ that is not in their complete class,
there is a design $\xi^*$ that is in the complete class such that all
elements of $I_{\xi^*}$ are the same as the corresponding elements
in~$I_{\xi}$, except that \textit{one diagonal element} in $I_{\xi
^*}$ is at\vadjust{\goodbreak}
least as large as that in $I_{\xi}$. This guarantees of course that
$I_{\xi^*} \geq I_{\xi}$. The contribution of this paper is that we
focus on increasing a \textit{principal submatrix} rather than just a
single \textit{diagonal element}.
This allows us to obtain results for more models than could be
addressed by \citet{ys09}, \citet{y09} and \citet{d05},
and also facilitates the identification of smaller complete
classes for some models considered in these earlier papers.

In Section~\ref{sec2} we will present the necessary background, while
the main results are featured in Section~\ref{sec3}. The power of the
proposed extension is seen through applications in Section~\ref{sec4}.
We conclude with a short discussion in Section~\ref{sec5}.\vspace*{-3pt}

\section{Information matrix and approximate designs}\label{sec2}
Consider a nonlinear regression model for which a response variable $y$
depends on a single regression variable $x$. We assume that the $y$'s
are independent and follow some exponential distribution $G$ with mean
$\eta(x,\theta)$, where $\theta$ is the $p\times1$ parameter vector,
and the values of $x$ can be chosen by the experimenter. Typically,
approximate designs are used to study optimality in this context. An
approximate design $\xi$ can be written as $\xi=\{(x_i,\omega_i),
i=1,\ldots,N\}$, where $\omega_i>0$ is the weight for design point
$x_i$ and $\sum_{i=1}^N \omega_i=1$. It is often more convenient to
present $\xi$ as $\xi=\{(c_i,\omega_i), i=1,\ldots,N\}$, $c_i\in
[A,B]$, with the $c_i$'s obtained from the $x_i$'s through a bijection
that may depend on $\theta$. Typically, the information matrix for
$\theta$ under design $\xi$ can be written as
\begin{equation}\label{infor1}
I_{\xi}(\theta)=P(\theta) \Biggl(\sum_{i=1}^N\omega_i C(\theta,
c_i) \Biggr)(P(\theta))^T,\vspace*{-1pt}
\end{equation}
where
\begin{equation}\label{infor2}
C(\theta, c)=
\pmatrix{
\Psi_{11}(c) & & & \cr
\Psi_{21}(c) & \Psi_{22}(c) & &\cr
\vdots& \vdots& \ddots& \cr
\Psi_{p1}(c) & \Psi_{p2}(c) & \cdots&\Psi_{pp}(c)
}.\vspace*{-1pt}
\end{equation}
The functions $\Psi$ are allowed to depend on $\theta$ not just through
$c$, but in an attempt to simplify notation we write, for example,
$\Psi
_{11}(c)$ rather than $\Psi_{11}(\theta, c)$. In (\ref{infor2}),
$C(\theta, c)$ is a symmetric matrix, and $P(\theta)$ is a $p\times p$
nonsingular matrix that depends only on $\theta$. Some examples of
(\ref
{infor1}) and (\ref{infor2}) will be seen in Section~\ref{sec4}.

For some $p_1$, $1\leq p_1 <p$, we partition $C(\theta, c)$ as
\begin{equation}\label{infor_part}
C(\theta, c)=
\pmatrix{
C_{11}(c)& C_{21}^T(c)\cr
C_{21}(c) & C_{22}(c)
}
.\vspace*{-1pt}
\end{equation}
Here, $C_{22}(c)$ is the lower $p_1\times p_1$ principal submatrix of
$C(\theta,c)$, that is,
\begin{equation}\label{infor_c22}
C_{22}( c)=
\pmatrix{
\Psi_{p-p_1+1,p-p_1+1}(c) & \cdots& \Psi_{p-p_1+1,p}(c)\cr
\vdots& \ddots& \vdots\cr
\Psi_{p, p-p_1+1}(c) & \cdots&\Psi_{pp}(c)
}
.\vadjust{\goodbreak}
\end{equation}

In the context of local optimality, if designs $\xi=\{(c_i,\omega_i),
i=1,\ldots,N\}$ and $\tilde{\xi}=\{(\tilde{c}_j,\tilde{\omega}_j),
j=1,\ldots,\tilde{N}\}$ satisfy $\sum_{i=1}^N\omega_i C(\theta,
c_i)\leq\sum_{i=1}^{\tilde{N}}\tilde{\omega}_i C(\theta, \tilde
{c}_i)$, then it follows from (\ref{infor1}) that $I_\xi(\theta) \le
I_{\tilde{\xi}}(\theta)$. Hence, $I_{\xi}(\theta)\leq I_{\tilde
{\xi
}}(\theta)$ follows if it holds that
\begin{eqnarray}\label{strategy:0}
\sum_{i=1}^N\omega_i C_{11}(c_i)&=& \sum_{i=1}^{\tilde{N}}\tilde
{\omega
}_i C_{11}( \tilde{c}_i), \nonumber\\
\sum_{i=1}^N\omega_i C_{12}( c_i)&=& \sum_{i=1}^{\tilde{N}}\tilde
{\omega
}_i C_{12}( \tilde{c}_i)\quad\mbox{and}\\
\sum_{i=1}^N\omega_i C_{22}( c_i)&\leq&\sum_{i=1}^{\tilde
{N}}\tilde
{\omega}_i C_{22}(\tilde{c}_i).\nonumber
\end{eqnarray}
This is what we explore in this paper. Note that this is more general
than \citet{ys09}, \citet{y09} and \citet{d05},
where $p_1=1$. We develop a theoretical framework for general values of $p_1$.

\section{Main results}\label{sec3}
Following \citet{k66} and \citet{d05}, a set
of $k+1$ real-valued continuous functions $u_0,\ldots,u_k$ defined on
an interval $[A,B]$ is called a Chebyshev system on $[A,B]$ if
\begin{equation}\label{cheby1}
\left|\matrix{
u_0(z_0) & u_0(z_1) & \cdots& u_0(z_k) \cr
u_1(z_0) & u_1(z_1) & \cdots& u_1(z_k) \cr
\vdots& \vdots& \ddots& \vdots\cr
u_{k}(z_0) & u_{k}(z_1) & \cdots& u_{k}(z_k)
}
\right|
\end{equation}
is strictly positive whenever $A\leq z_0<z_1<\cdots<z_k\leq B$.

Along the lines of \citet{y09}, we select a maximal set of linearly
independent nonconstant functions from the $\Psi$ functions that appear
in the first $p-p_1$ columns of the matrix $C(\theta,c)$ defined in
(\ref{infor2}), and rename the selected functions as $\Psi_1, \ldots,
\Psi_{k-1}$. For a given nonzero $p_1 \times1$ vector $Q$, let
\begin{equation}\label{PsiQ}
\Psi^Q_{k}=Q^TC_{22}(c)Q,
\end{equation}
where $C_{22}(c)$ is as defined in (\ref{infor_c22}).

For $\Psi_0 = 1$, $\Psi_1, \ldots, \Psi_{k-1}$ and $C_{22}(c)$, we will
say that a set of $n_1$ pairs $(c_i, \omega_i)$ is dominated by a set
of $n_2$ pairs $(\tilde{c}_i, \tilde{\omega}_i)$ if
\begin{eqnarray}\label{exist:1}\quad
\sum_i\omega_i\Psi_l(c_{i})&=&\sum_i\tilde{\omega}_i\Psi
_l(\tilde{c}_{i}),\qquad
l=0,1,\ldots,k-1; \\
\label{exist:2}
\sum_i\omega_i\Psi^Q_{k}(c_{i})&<&\sum_i\tilde{\omega}_i\Psi
^Q_{k}(\tilde
{c}_{i})\qquad\mbox{for every nonzero vector $Q$},
\end{eqnarray}
where the summations on the left-hand sides are over the $n_1$
subscripts for the pairs $(c_i, \omega_i)$ and those on the right-hand
sides over the $n_2$ subscripts for the pairs $(\tilde{c}_i, \tilde
{\omega}_i)$.

The following two lemmas provide the basic tools for the main results.
We point out that the pairs $(c_i, \omega_i)$ in these lemmas need not
form a design; in particular, the $\omega_i$'s need not add to 1.

\begin{Lemma} \label{lemma1}
For the functions $\Psi_0=1, \Psi_1, \ldots, \Psi_{k-1}, \Psi^Q_k$
defined on an interval $[A,B]$, suppose that either
\begin{eqnarray}\label{assumption1}
&&\{\Psi_0, \Psi_1, \ldots, \Psi_{k-1}\} \mbox{ and } \{\Psi_0,
\Psi_1,
\ldots, \Psi_{k-1}, \Psi^Q_k\} \nonumber\\[-8pt]\\[-8pt]
&&\mbox{form Chebyshev systems for every nonzero vector } Q\nonumber
\end{eqnarray}
or
\begin{eqnarray}\label{assumption2}
&&\{\Psi_0, \Psi_1, \ldots, \Psi_{k-1}\} \mbox{ and } \{\Psi_0,
\Psi_1,
\ldots, \Psi_{k-1}, -\Psi^Q_k\} \nonumber\\[-8pt]\\[-8pt]
&&\mbox{form Chebyshev systems for every nonzero vector } Q.\nonumber
\end{eqnarray}

Then the following conclusions hold:
\begin{longlist}[(b)]
\item[(a)] For $k=2n-1$, if (\ref{assumption1}) holds, then for any
set $S_1 = \{(c_i, \omega_i)\dvtx  \omega_i>0,\break i=1,\ldots,n$\} with
$A\leq
c_1<\cdots<c_n< B$, there exists a set $S_2 = \{(\tilde{c}_i,\tilde
{\omega}_i)\dvtx\break \tilde{\omega}_i>0, i=1,\ldots,n\}$ with $c_1<\tilde
{c}_1<c_2 <\cdots<\tilde{c}_{n-1}<c_n<\tilde{c}_n=B$, such that~$S_1$
is dominated by $S_2$.

\item[(b)] For $k=2n-1$, if (\ref{assumption2}) holds, then for any
set $S_1 = \{(c_i, \omega_i)\dvtx \omega_i>0,\break i=1,\ldots,n\}$ with
$A<c_1<\cdots<c_n\leq B$, there exists a set $S_2 = \{(\tilde
{c}_i,\tilde{\omega}_i)\dvtx \break \tilde{\omega}_i>0, i=0,\ldots,n-1\}$ with
$A= \tilde{c}_0<c_1<\tilde{c}_1<c_2<\cdots<\tilde{c}_{n-1}<c_n$, such
that $S_1$ is dominated by $S_2$.

\item[(c)] For $k=2n$, if (\ref{assumption1}) holds, then for any set
$S_1 = \{(c_i, \omega_i)\dvtx \omega_i>0, \break i=1,\ldots,n\}$ with $A <
c_1<\cdots<c_n< B$, there exists a set $S_2 = \{(\tilde{c}_i,\tilde
{\omega}_i)\dvtx\break  \tilde{\omega}_i>0,  i=0,\ldots,n\}$ with $ A=\tilde
{c}_0<c_1<\tilde{c}_1< \cdots<c_n<\tilde{c}_n=B$, such that $S_1$ is
dominated by $S_2$.

\item[(d)] For $k=2n$, if (\ref{assumption2}) holds, then for any set
$S_1 = \{(c_i, \omega_i), \omega_i>0,\break i=1,\ldots,n+1$ with $A\leq
c_1<\cdots<c_{n+1} \leq B$, there exists a set $S_2 = \{(\tilde
{c}_i,\tilde{\omega}_i)\dvtx\break \tilde{\omega}_i>0,  i=1,\ldots,n\}$ with
$c_1<\tilde{c}_1< \cdots<c_n<\tilde{c}_n <c_{n+1}$, such that $S_1$ is
dominated by~$S_2$.
\end{longlist}
\end{Lemma}

\begin{pf}
Since the proof is similar for all parts, we only provide a proof for
part (a).

Let $S_1$ be as in part (a). First consider the special case that
$Q\!=\!(1,0,\ldots,0)^T$. By (1a) of Therorem 3.1 in \citet{d05}, there exists a set of at most $n$ pairs ($\tilde{c}_i, \tilde
{\omega}_i)$ with one of the points equal to $B$ so that~(\ref
{exist:1}) and~(\ref{exist:2}) hold for this $Q$. By part (a) of
Proposition~\ref{prop0} in the \hyperref[append]{Appendix}, the number of distinct points with
$\tilde{\omega}_i >0$ must then be exactly $n$. Thus we have\vadjust{\goodbreak} $\tilde
{c}_1 < \cdots< \tilde{c}_n=B$, and the $c_i$'s and $\tilde{c}_i$'s
must alternate by part (b) of Proposition~\ref{prop0}. The result follows now
for an arbitrary nonzero $Q$ by applying Proposition~\ref{prop1} in the
\hyperref[append]{Appendix}
and using (\ref{assumption1}) and (\ref{exist:2}).
\end{pf}

Lemma~\ref{lemma2} partially extends Lemma~\ref{lemma1} by observing
that larger sets $S_1$ than in Lemma~\ref{lemma1} are also dominated by
sets $S_2$ as in that lemma.

\begin{Lemma} \label{lemma2}
With the same notation and assumptions as in Lemma~\ref{lemma1}, let
$S_1 = \{(c_i,\omega_i)\dvtx \omega_i>0, A\le c_i \le B,  i=1,\ldots,N\}$,
where $N \geq n$ for cases (\textup{a}), (\textup{b}), and (\textup{c}) of Lemma~\ref
{lemma1}, and $N \geq n+1$ for case (\textup{d}). Then the following
conclusions hold:
\begin{longlist}[(b)]
\item[(a)] For $k=2n-1$, if (\ref{assumption1}) holds, then $S_1$ is
dominated by a set $S_2$ of size $n$ that includes $B$ as one of the points.
\item[(b)] For $k=2n-1$, if (\ref{assumption2}) holds, then $S_1$ is
dominated by a set $S_2$ of size $n$ that includes $A$ as one of the points.
\item[(c)] For $k=2n$, if (\ref{assumption1}) holds, then $S_1$ is
dominated by a set $S_2$ of size $n+1$ that includes both $A$ and $B$
as points.
\item[(d)] For $k=2n$, if (\ref{assumption2}) holds, then $S_1$ is
dominated by a set $S_2$ of size $n$.
\end{longlist}
\end{Lemma}

\begin{pf}
The results follow by application of Lemma~\ref{lemma1}. For example,
for case~(a), if $N=n$, the result follows directly from Lemma~\ref
{lemma1}. If $N>n$, we start with the points $c_1 < c_2 < \cdots< c_N$
in $S_1$. Using Lemma~\ref{lemma1}, we obtain points $c_1, \ldots,
c_{N-n}, \tilde{c}_{N-n+1}, \ldots, \tilde{c}_N=B$ in a set $\tilde
{S}_1$ that dominates $S_1$. Using Lemma~\ref{lemma1} again on the $n$
largest points other than $\tilde{c}_N$ in $\tilde{S}_1$, we move one
more point to $B$, obtaining a new set with $N-1$ points that dominates
$\tilde{S}_1$. Continue until the size of the set is reduced to $n$;
this is the desired set $S_2$.
\end{pf}

The first main result is an immediate consequence of Lemma~\ref{lemma2}.

\begin{Theorem} \label{main1}
For a regression model with a single regression variable $x$, suppose
that the information matrix $C(\theta,c)$ can be written as in (\ref
{infor1}) for $c \in[A,B]$. Partitioning the information matrix as in
(\ref{infor_part}), let $\Psi_1, \ldots, \Psi_{k-1}$ be a maximum set
of linearly independent nonconstant $\Psi$ functions in the first
$p-p_1$ columns of $C(\theta, c)$. Define $\Psi^Q_{k}$ as in (\ref
{PsiQ}). Suppose that either~(\ref{assumption1}) or (\ref{assumption2})
in Lemma~\ref{lemma1} holds. Then the following complete class results hold:
\begin{longlist}[(b)]
\item[(a)] For $k=2n-1$, if (\ref{assumption1}) holds, the designs
with at most $n$ support points, including $B$, form a complete class.
\item[(b)] For $k=2n-1$, if (\ref{assumption2}) holds, the designs
with at most $n$ support points, including $A$, form a complete class.
\item[(c)] For $k=2n$, if (\ref{assumption1}) holds, the designs with
at most $n+1$ support points, including both $A$ and $B$, form a
complete class.
\item[(d)] For $k=2n$, if (\ref{assumption2}) holds, the designs with
at most $n$ support points form a complete class.\vadjust{\goodbreak}
\end{longlist}
\end{Theorem}

Note that if (\ref{exist:1}) holds for $\Psi_l(c)$, $l = 1,\ldots,k-1$,
then the same is true if we replace one or more of the $\Psi_l$'s by
$-\Psi_l$. Therefore, if (\ref{assumption1}) or
(\ref{assumption2}) do not hold for the original $\Psi_l$'s,
conclusions in Theorem~\ref{main1} would still be valid if~(\ref
{assumption1}) and (\ref{assumption2}) hold after multiplying one or
more of the $\Psi_l$'s, $l=1,\ldots,k-1$, by $-1$.

While Theorem~\ref{main1} is very powerful, applying it directly may
not be easy. The next result, which utilizes a generalization of a tool
in \citet{y09}, will lead to a condition that is easier to verify.
Using the notation of Theorem~\ref{main1}, define functions $f_{l,t}$,
$1\leq t\leq k; t\leq l\leq k$ as follows:
\begin{equation}\label{def:df}
f_{l,t}(c)= \cases{
\Psi_l'(c), & \mbox{if $t=1$, $l=1,\ldots,k-1$}, \cr
C'_{22}(c), & \mbox{if $t=1$, $l=k$}, \cr
\biggl(\dfrac{f_{l,t-1}(c)}{f_{t-1,t-1}(c)}\biggr )', & \mbox{if $2\leq
t\leq k$, $t\leq l\leq k$.}
}
\end{equation}
The following lower triangular matrix contains all of these functions,
and suggest an order in which to compute them:
\begin{equation}\label{df}
\pmatrix{
f_{1,1}=\Psi_1' & & & & \cr
f_{2,1}=\Psi_2' & f_{2,2}= \big(\frac{f_{2,1}}{f_{1,1}}\big )' & & &
\vspace*{2pt}\cr
f_{3,1}=\Psi_3' & f_{3,2}= \big(\frac{f_{3,1}}{f_{1,1}} \big)' &
f_{3,3}= \big(\frac{f_{3,2}}{f_{2,2}} \big)' & &
\cr
\vdots& \vdots& \vdots& \ddots&\cr
f_{k,1}=C'_{22} & f_{k,2}= \big(\frac{f_{k,1}}{f_{1,1}} \big)' &
f_{k,3}= \big(\frac{f_{k,2}}{f_{2,2}} \big)' & \vdots& f_{k,k}=
\big(\frac{f_{k,k-1}}{f_{k-1,k-1}} \big)'
}
.\hspace*{-35pt}
\end{equation}
Note that, for $p_1 \ge2$, the functions in the last row are matrix
functions, which is a key difference with \citet{y09}. The derivatives
of matrices in~(\ref{def:df}) are element-wise derivatives. For the
next result, we will make the following assumptions:
\begin{longlist}[(ii)]
\item[(i)] All functions $\Psi$ in the information matrix $C (\theta
,c)$ are at least $k$th order differentiable on $(A,B)$.
\item[(ii)] For $1\leq l\leq k-1$, the functions $f_{l,l}(c)$ have no
roots in $[A,B]$.
\end{longlist}

For ease of notation, in the remainder we will write $f_{l,l}$ instead
of $f_{l,l}(c)$, and $f_{l,l}>0$ means that $f_{l,l}(c)>0$ for all
$c\in
[A,B]$. This also applies for $l=k$, in which case it means that the
matrix $f_{k,k}$ is positive definite for all $c \in[A,B]$.

\begin{Theorem} \label{main2}
For a regression model with a single regression variable~$x$, let $c
\in[A,B]$, $C(\theta,c)$, $\Psi_1,\ldots,\Psi_{k-1}$ and $\Psi
_k^Q$ be
as in Theorem~\ref{main1}. For the functions $f_{l,l}$ in (\ref
{def:df}), define $F(c)=\prod_{l=1}^kf_{l,l}$, $c\in[A,B]$. Suppose
that either $F(c)$ or $-F(c)$ is positive definite for all $c\in
[A,B]$. Then the following complete class results hold:
\begin{longlist}[(b)]
\item[(a)] For $k=2n-1$, if $F(c)>0$, the designs with at most $n$
support points, including $B$, form a complete class.
\item[(b)] For $k=2n-1$, if $-F(c)>0$, the designs with at most $n$
support points, including $A$, form a complete class.
\item[(c)] For $k=2n$, if $F(c)>0$, the designs with at most $n+1$
support points, including both $A$ and $B$, form a complete class.
\item[(d)] For $k=2n$, if $-F(c)>0$, the designs with at most $n$
support points form a complete class.
\end{longlist}
\end{Theorem}

\begin{pf}
We only present the proof for case (a) since the other cases are
similar. For any nonzero vector $Q$, $Q^T F(c) Q>0$ for all $c\in
[A,B]$. Among all $f_{l,l}$, $l=1,\ldots,k-1$, and $Q^T f_{k,k}Q$,
suppose that $a$ of them are negative. Let $1 \le l_1 < \cdots< l_a
\le k$ denote the subscripts for these negative terms, and note that
$a$ must be even. Note also that the labels $l_1 < \cdots< l_a$ do not
depend on the choice of the vector $Q$ since $f_{1,1}, \ldots,
f_{k-1,k-1}$ do not depend on $Q$. Finally, note that for any $l$ with
$1 \le l \le k-1$, if we replace $\Psi_l(c)$ by $-\Psi_l(c)$, then the
signs of $f_{l,l}$ and $f_{l+1,l+1}$ are switched while all others
remain unchanged.

We now change some of the $\Psi_l$'s to $-\Psi_l$. This is done for
those $l$ that satisfy $l_{2b-1}\leq l <l_{2b}$ for some value of $b
\in\{1, \ldots, a/2\}$. Denote the new $\Psi$-functions by $\{
1,\widehat{\Psi}_1,\ldots, \widehat{\Psi}_k^Q \}$. Notice that
$\widehat
{\Psi}_k^Q=\Psi_k^Q$. From the last observation in the previous
paragraph, it is easy to check that $f_{l,l}>0$, $l=1,\ldots,k$, for
the functions $f_{l,l}$ that correspond to this new set of $\widehat
{\Psi}$-functions. By Proposition~\ref{prop3} in the \hyperref[append]{Appendix}, $\{1,
\widehat{\Psi}_1, \ldots, \widehat{\Psi}_{k-1}\}$ and $\{1,
\widehat
{\Psi}_1, \ldots, \widehat{\Psi}_{k-1}, \widehat{\Psi}_{k}^Q\}$ are
Chebyshev systems on $[A,B]$, regardless of the choice for $Q \not=
0$. The result follows now from case (a) of Theorem~\ref{main1} and
the observation immediately after Theorem~\ref{main1}.
\end{pf}

For case (a) in Theorem~\ref{main2}, the value of $A$ in the interval
$[A,B]$ is allowed to be $-\infty$. In this situation, for any given
design $\xi$, we can choose $A=\min_{i} c_i$, and the conclusion of the
theorem holds. Similarly, $B$ can be $\infty$ in case (b), and the
interval can be unbounded at either side for case (d).

As noted at the end of Section~\ref{sec2}, the results in \citet{ys09},
\citet{y09}, and \citet{d05} correspond to
$p_1=1$. The extension in this paper allows the choice of larger values
of $p_1$ where feasible. Larger values of $p_1$ lead to designs with
smaller support sizes. The reason for this is that the value of $k$ in
Theorems~\ref{main1} and~\ref{main2} corresponds to the number of
equations in (\ref{exist:1}). For a particular model, this number is
smaller for larger~$p_1$. Since the support size of the designs is
roughly half the value of~$k$, the support size is smaller for larger
values of $p_1$.

We will provide some examples of the application of Theorems~\ref
{main1} and~\ref{main2} in the next section, and will offer some
further thoughts on the ease of their application in Section~\ref{sec5}.

\section{Applications}\label{sec4}
Whether the model is for continuous or discrete data, with homogeneous
or heterogeneous errors, Theorems~\ref{main1} and~\ref{main2} can be
applied as long as the information matrix can be written as in (\ref
{infor1}). As the examples in this section will show, in many cases the
result of the theorem facilitates the determination of complete classes
with the minimal number of support points.

\subsection{Exponential regression models}

Dette, Melas and Wong (\citeyear{dmw06}) studied exponential regression models,
which can be written as
\begin{equation}\label{example1}
Y_i=\sum_{l=1}^L a_l e^{-\lambda_l x_i}+\varepsilon_i,
\end{equation}
where the $\varepsilon_{i}$'s are i.i.d. with mean 0 and variance
$\sigma
^2$, and $x_i \in[U,V]$ is the value of the regression variable to be
selected by the experimenter. Here
$\theta=(a_1,\ldots,a_L,\lambda_1,\ldots,\lambda_L)^T$, with
$a_l\ne
0$, $l=1,\ldots,L$, and
$0<\lambda_1<\cdots<\lambda_L$. For $L=2$, they showed that there is a
$D$-optimal design for
$\theta=(a_1, a_2,\lambda_1,\lambda_2)^T$ based on four points,
including the lower limit $U$. Further, for $L=3$ and $\lambda
_2=(\lambda_1+\lambda_3)/2$, they showed that there is a $D$-optimal
design for $\theta$ based on six points, again including the lower
limit $U$. By using Theorem~\ref{main2}, we will show that similar
conclusions are possible for other optimality criteria, including $A$-
and $E$-optimality, and other functions of interest for many a~priori
values of $\theta$.

For $L=2$, the results in \citet{y09} can be used to obtain a complete
class of designs with at most five points. We can do better with
Theorem~\ref{main2}. The information matrix for $\theta
=(a_1,a_2,\lambda
_1,\lambda_2)^T$
under design $\{(x_i,\omega_i), i=1,\ldots,N\}$ can be written in the
form of (\ref{infor1}) with $P(\theta)=\operatorname{diag}(1,1,\frac{a_1}{\lambda
_2-\lambda_1},\frac{a_2}{\lambda_2-\lambda_1})$ and
\begin{equation}\label{example_1}
\quad C(\theta,c)=
\pmatrix{
c^\lambda& & & \cr
c^{\lambda+1} & c^{\lambda+2} & & \cr
\log(c)c^{\lambda}& \log(c)c^{\lambda+1} & \log^2(c)c^{\lambda} & \vspace*{2pt}\cr
\log(c)c^{\lambda+1} & \log(c)c^{\lambda+2} & \log^2(c)c^{\lambda+1} &
\log^2(c)c^{\lambda+2}
}
,
\end{equation}
where $c=e^{-(\lambda_2-\lambda_1)x}$ and $\lambda=\frac{2\lambda
_1}{\lambda_2-\lambda_1}$. Let $\Psi_1(c)=c^\lambda$, $\Psi
_2(c)=\log(c)c^{\lambda}$,\break $\Psi_3(c)=c^{\lambda+1}$, $\Psi
_4(c)=\log(c)c^{\lambda+1}$, $\Psi_5(c)=c^{\lambda+2}$, $\Psi
_6(c)=\log(c)c^{\lambda+2}$ and
\[
C_{22}(c)=
\pmatrix{
\log^2(c)c^{\lambda} & \log^2(c)c^{\lambda+1}\vspace*{2pt}\cr
\log^2(c)c^{\lambda+1} & \log^2(c)c^{\lambda+2}
}
.
\]
Then $f_{1,1}=\lambda c^{\lambda-1}$, $f_{2,2}=\frac{1}{c}$,
$f_{3,3}=\frac{\lambda+1}{\lambda}$, $f_{4,4}=\frac{1}{c}$,
$f_{5,5}=\frac{4(\lambda+2)}{\lambda+1}$, $f_{6,6}=\frac{1}{c}$ and
\[
f_{7,7}(c)=
\pmatrix{
\dfrac{2\lambda}{(\lambda+2)c^3} & \dfrac{\lambda+1}{2(\lambda
+2)c^2}\vspace*{2pt}\cr
\dfrac{\lambda+1}{2(\lambda+2)c^2} & \dfrac{2}{c}
}
.
\]
Note that $c>0$ and $\lambda>0$, so that $F(c)$ is positive definite if
$|f_{7,7}(c)|>0$. This is equivalent to $15\lambda^2+30\lambda-1>0$,
which is satisfied when $\frac{\lambda_2}{\lambda_1}<\frac{\sqrt
{960}+30}{\sqrt{960}-30}$. Thus, by (a) of Theorem~\ref{main2}, we
have the following result.

\begin{Theorem}
For Model (\ref{example1}) with $L=2$, if
\[
\frac{\lambda_2}{\lambda_1}<\frac{\sqrt{960}+30}{\sqrt
{960}-30}\thickapprox
61.98,
\]
then the designs with at most four points, including the lower limit
$U$, form a complete class.
\end{Theorem}

For $L=3$ and $2\lambda_2=\lambda_1+\lambda_3$, the information matrix
for $\theta=(a_1,a_2,a_3,\break\lambda_1,\lambda_2,\lambda_3)^T$ under design
$\{(x_i,\omega_i), i=1,\ldots,N\}$ can be written in the form of
(\ref{infor1}) with $P(\theta)=\operatorname{diag}(1,1,1, \frac{a_1}{\lambda
_2-\lambda
_1},\frac{a_2}{\lambda_2-\lambda_1},\frac{a_3}{\lambda_2-\lambda_1})$
and
\begin{equation}\label{example12}
 C(\theta,c)\!=\!
{\fontsize{7.8}{10}\selectfont{
\pmatrix{
c^\lambda& \hspace*{-2.5pt}&\hspace*{-2.5pt} &\hspace*{-2.5pt} &\hspace*{-2.5pt} & \cr
c^{\lambda+1} &\hspace*{-2.5pt} c^{\lambda+2} &\hspace*{-2.5pt} &\hspace*{-2.5pt} &\hspace*{-2.5pt} & \cr
c^{\lambda+2} &\hspace*{-2.5pt} c^{\lambda+3} &\hspace*{-2.5pt} c^{\lambda+4} &\hspace*{-2.5pt} &\hspace*{-2.5pt} & \cr
\log(c)c^{\lambda} &\hspace*{-2.5pt} \log(c)c^{\lambda+1} &\hspace*{-2.5pt} \log(c)c^{\lambda+2} &\hspace*{-2.5pt} \log^2(c)c^{\lambda} &\hspace*{-2.5pt} & \vspace*{2pt}\cr
\log(c)c^{\lambda+1} &\hspace*{-2.5pt} \log(c)c^{\lambda+2} &\hspace*{-2.5pt} \log(c)c^{\lambda+3} &\hspace*{-2.5pt} \log^2(c)c^{\lambda+1} &\hspace*{-2.5pt} \log^2(c)c^{\lambda+2}& \vspace*{2pt}\cr
\log(c)c^{\lambda+2} &\hspace*{-2.5pt} \log(c)c^{\lambda+3} &\hspace*{-2.5pt} \log(c)c^{\lambda+4} &\hspace*{-2.5pt}\log^2(c)c^{\lambda+2} &\hspace*{-2.5pt} \log^2(c)c^{\lambda+3}&\hspace*{-2.5pt} \log^2(c)c^{\lambda+4}}}},
\hspace*{-35pt}
\end{equation}
where $c=e^{-(\lambda_2-\lambda_1)x}$ and $\lambda=\frac{2\lambda
_1}{\lambda_2-\lambda_1}$. Let $\Psi_{2l-1}(c)=c^{\lambda+l-1}$ and
$\Psi_{2l}(c)=\log(c)c^{\lambda+l-1}$, $l=1,\ldots,5$, and let
\[
C_{22}(c)=
\pmatrix{
\log^2(c)c^{\lambda} & \log^2(c)c^{\lambda+1} & \log^2(c)c^{\lambda+2}
\vspace*{2pt}\cr
\log^2(c)c^{\lambda+1} & \log^2(c)c^{\lambda+2}& \log^2(c)c^{\lambda
+3}\vspace*{2pt}\cr
\log^2(c)c^{\lambda+2} & \log^2(c)c^{\lambda+3}& \log^2(c)c^{\lambda+4}
}
.
\]
Then $f_{1,1}=\lambda c^{\lambda-1}$, $f_{2l,2l}=\frac{1}{c}$,
$l=1,2,3,4,5$, $f_{2l+1,2l+1}=\frac{l^2(\lambda+l)}{\lambda+l-1}$,
$l=1,2,3,4$, and
\[
f_{11,11}(c)=
\pmatrix{
\dfrac{2\lambda}{(\lambda+4)c^5} & \dfrac{\lambda+1}{8(\lambda
+4)c^4} &
\dfrac{\lambda+2}{18(\lambda+4)c^3}\vspace*{2pt}\cr
\dfrac{\lambda+1}{8(\lambda+4)c^4} & \dfrac{\lambda+2}{18(\lambda+4)c^3}
& \dfrac{\lambda+3}{8(\lambda+4)c^2}\vspace*{2pt}\cr
\dfrac{\lambda+2}{18(\lambda+4)c^3} &
\dfrac{\lambda+3}{8(\lambda+4)c^2} &\dfrac{2}{c}
}
.
\]
Again, $c>0$ and $\lambda>0$, so that $F(c)$ is positive definite if
$|(f_{11,11}(c))|$ and its leading principal minors are positive. This
is equivalent to
\begin{eqnarray}\label{example12:1}
1505\lambda^3+9030\lambda^2+11499\lambda-1082&>&0,\nonumber\\
55\lambda^2+110\lambda-9&>&0,\nonumber\\[-8pt]\\[-8pt]
1295\lambda^2+5180\lambda- 4&>&0\nonumber\\
\mbox{and}\quad 55\lambda^2 + 330\lambda+ 431&>&0.\nonumber
\end{eqnarray}
Simple computation shows that this holds for $\frac{\lambda
_2}{\lambda
_1}<23.72$ (or, equivalently, $\frac{\lambda_3}{\lambda_1}<46.45$). By
Theorem~\ref{main2}, we have the following result.

\begin{Theorem}
For model~(\ref{example1}) with $L=3$ and $2\lambda_2=\lambda
_1+\lambda
_3$, if $\frac{\lambda_2}{\lambda_1}<23.72$, then the designs with at
most six points, including the lower limit $U$, form a complete class.
\end{Theorem}

\subsection{LINEXP model}
Demidenko (\citeyear{De2006}) proposed a model referred to as the LINEXP model to
describe tumor growth delay and regrowth. The natural logarithm of the
tumor volume is modeled as
\begin{equation}\label{linexp}
Y_{i}=\alpha+\gamma x_i +\beta(e^{-\delta x_i}-1)+\varepsilon_{i},
\end{equation}
with independent $\varepsilon_{i}\sim N(0,\sigma^2)$ and $x_i\in[U, V]$
as the value of the single regression variable, which in this case
refers to time. Here $\theta=(\alpha, \gamma, \beta, \delta)^T$ is the
parameter vector, where $\alpha$ is the baseline logarithm of the tumor
volume, $\gamma$ is the final growth rate and $\delta$ is the rate at
which killed cells get washed out. The size of the parameter $\beta$
relative to $\gamma/\delta$ determines whether regrowth is monotonic
($\beta< \gamma/\delta$) or not. \citet{l08} recently
studied this model and showed that a $D$-optimal design for $\theta$
can be based on four points, including $U$ and $V$. We will now show
that Theorem~\ref{main2} extends this conclusion to other optimality
criteria and functions of interest.

The information matrix for $\theta$ under design $\{(x_i,\omega_i),
i=1,\ldots,N\}$ can be written in the form of (\ref{infor1}) with
\begin{eqnarray}\label{example2}
P(\theta)&=&
\pmatrix{
1 & 0 & 0 & 0\cr
1 & 0 &1 & 0 \cr
0 & -\delta& 0 & 0\cr
0 & 0 & 0 & \delta/\beta
}
^{-1}\quad
\mbox{and}\nonumber\\[-8pt]\\[-8pt]
C(\theta,c)&=&
\pmatrix{
1 & & & \cr
e^{c} & e^{2c} & & \cr
c & c e^{c} &c^2 & \cr
c e^{c} & c e^{2c} & c^2 e^{c} & c^2 e^{2c}
}
,\nonumber
\end{eqnarray}
where $c=-\delta x$. With a proper choice of $\Psi$ functions, it can
be shown that the result in \citet{y09} yields a complete class of
designs with at most five points, including $U$ and $V$. We can again
do better with Theorem~\ref{main2}.

Define $\Psi_1(c)=c$, $\Psi_2(c)=e^c$, $\Psi_3(c)=ce^c$, $\Psi
_4(c)=e^{2c}$, $\Psi_5(c)=ce^{2c}$ and
\[
C_{22}(c)=
\pmatrix{
c^2 & c^2 e^c\cr
c^2 e^c & c^2 e^{2c}
}
.
\]
This yields $f_{1,1}=1$, $f_{2,2}=e^c$, $f_{3,3}=1$, $f_{4,4}=4e^c$,
$f_{5,5}=1$ and
\[
f_{6,6}(c)=
\pmatrix{
2e^{-2c}&e^{-c}/2\cr
e^{-c}/2 & 2
}
.
\]
Clearly $F(c)$ is a positive definite matrix. Therefore, by part (c)
of Theorem~\ref{main2}, we reach the following conclusion.

\begin{Theorem}
For the LINEXP model (\ref{linexp}), the designs with at most four
points, including $U$ and $V$, form a complete class.
\end{Theorem}

\subsection{Double-exponential regrowth model} \label{subsec3}

Demidenko (\citeyear{d04}), using a two-compartment model, developed a
double-exponential regrowth model to describe the dynamics of
post-irradiated tumors. The model can be written as
\begin{equation}\label{double}
Y_{i}=\alpha+\ln[\beta e^{\nu x_i}+(1-\beta)e^{-\phi x_i}]
+\varepsilon_{ij},
\end{equation}
with independent $\varepsilon_{i}\sim N(0,\sigma^2)$ and $x_i \in[U, V]$
again as the value for the variable time. Here $\theta=(\alpha, \beta
,\nu,\phi)^T$ is the parameter vector, where $\alpha$ is the logarithm
of the initial tumor volume, $0<\beta<1$ is the proportional
contribution of the first compartment and $\nu$ and $\phi$ are cell
proliferation and death rates.

Using Chebyshev systems and an equivalence theorem, \citet{l08}
showed that a $D$-optimal design for $\theta$ can be based on
four points including $U$ and $V$. Theorem~\ref{main1} allows us to
extend this result to a~complete class result, thereby covering many
other optimality criteria and any functions of interest.

The information matrix for $\theta$ under design $\{(x_i,\omega_i),
i=1,\ldots,N\}$ is of the form (\ref{infor1}) with
\[
P(\theta)=
\pmatrix{
1 & 0 & 0 & 0\cr
1 & 1-\beta& 0 & 0\cr
0 & 0 &1/\beta& 0 \cr
0 & 0 & 0 & -1/(1-\beta)
}
^{-1}
\]
and with $C(\theta,x)$ a $4\times4$ matrix as in (\ref{infor2}), where
$\Psi_{11}=1$, $\Psi_{21}=e^{\nu x}/g(x)$, $\Psi_{22}=e^{2\nu
x}/g^2(x)$, $\Psi_{31}=x e^{\nu x}/g(x)$, $\Psi_{32}=x e^{2\nu
x}/g^2(x)$, $\Psi_{33}=x^2e^{2\nu x}/g^2(x)$, $\Psi_{41}=x e^{-\phi
x}/g(x)$, $\Psi_{42}=x e^{(\nu-\phi)x}/g^2(x)$, $\Psi
_{43}=x^2e^{(\nu
-\phi)x}/g^2(x)$ and $\Psi_{44}=x^2e^{-2\phi x}/g^2(x)$. Here,
$g(x)=\beta e^{\nu x}+(1-\beta)e^{-\phi x}$. Note that $\Psi_{42}$ can
be written as a linear combination of $\Psi_{31}$ and $\Psi_{32}$. We
can apply Theorem~\ref{main1} if we can show that both $\{1, \Psi_{21},
\Psi_{22}, \Psi_{41}, -\Psi_{31}, \Psi_{32}\}$ and $\{1, \Psi
_{21}, \Psi
_{22}, \Psi_{41}$, $-\Psi_{31}, \Psi_{32}, Q^TC_{22}(x)Q\}$ are
Chebyshev systems for any nonzero vector $Q$, where $C_{22}(x)=
\left({
\Psi_{33} \atop \Psi_{43}}\enskip{
\Psi_{43} \atop \Psi_{44}}
\right)
$.

Rather than do this directly, we first simplify the problem. We
multiply each of the $\Psi$'s by the positive function $e^{2\phi
x}g(x)^2$, which preserves the Chebyshev system property. After further
simplifications by replacing some of the resulting functions by
independent linear combinations of these functions, which also
preserves the Chebyshev system property, we arrive at the systems\vadjust{\goodbreak} $\{
1$, $e^{(\nu+\phi)x}$, $e^{2(\nu+\phi)x}$, $x$, $-xe^{(\nu+\phi)x}$,
$xe^{2(\nu+\phi)x}\}$ and $\{1$, $e^{(\nu+\phi)x}$, $e^{2(\nu+\phi)x}$,
$x$, $-xe^{(\nu+\phi)x}$, $xe^{2(\nu+\phi)x}$, $g^2(x)e^{2\phi
x}Q^TC_{22}(x)Q\}$. It suffices to show that these are Chebyshev
systems for any nonzero vector $Q$, which follows from Proposition~\ref
{prop3} if we show that $f_{l,l}>0$, $l=1,\ldots,6$, for the latter
system. It can be shown that $f_{1,1}=f_{2,2}/2=2f_{4,4}=f_{5,5}/4=a
e^{ax}$, $f_{3,3}=e^{-2ax}$ and $f_{6,6}=
\left({2 \atop e^{-ax}/2}\enskip{
e^{-ax}/2 \atop 2e^{-2ax}}
\right)
$, where $a=\nu+\phi$. Thus both systems are Chebyshev systems, and by
part (c) of Theorem~\ref{main1}, we reach the following
conclusion.\vspace*{-2pt}

\begin{Theorem}
For the double-exponential regrowth model (\ref{double}), the designs
with at most four points, including $U$ and $V$, form a complete class.\vspace*{-2pt}
\end{Theorem}

\section{Discussion}\label{sec5}
We have given a powerful extension of the result in \citet{y09} that
has potential for providing a small complete class of designs whenever
the information matrix can be written as in (\ref{infor1}).
Irrespective of the optimality criterion (provided that it does not
violate the Loewner ordering) and of the function of $\theta$ that is
of interest, the search for an optimal design can be restricted to the
small complete class. As the examples in Section~\ref{sec4} show, the
results lead us to conclusions that were not possible using the results
in \citet{y09} and \citet{d05}.

As already pointed out, direct application of Theorem~\ref{main1} may
not be easy. Section~\ref{subsec3} shows some tricks that can be
useful when using Theorem~\ref{main1}. Direct application of
Theorem~\ref{main2} is easier because the condition for the
function~$F(c)$ can be verified with the help of software for symbolic
computations. Sometimes it is more convenient to do this after
multiplying each of the $\Psi$ functions by the same positive function
(see Section~\ref{subsec3}).

There remain, however, some basic questions related the application of
either Theorem~\ref{main1} or Theorem~\ref{main2} that do not have
simple general answers. For example, what is a good choice for $p_1$ in
forming the matrix $C_{22}(c)$ in (\ref{infor_c22})? In Section~\ref
{sec4}, the choice $p_1 = p/2$ worked well, and selecting $p_1$
approximately equal to $p/2$ may be a good general starting point.
Moreover, there is the question of how to order the rows and columns of
the information matrix. By reordering the elements in the parameter
vector $\theta$, we could wind up with different matrices $C_{22}(c)$,
even after fixing $p_1$. So what ordering is best? In all of the
examples in Section~\ref{sec4}, we have used an ordering that makes
``higher-order terms'' appear in $C_{22}(c)$, and this may offer the
best general strategy. There is still another issue related to
ordering: In renaming the independent $\Psi$-functions in the first
$p-p_1$ columns of $C(\theta,c)$, different orders will result in
different $f_{l,l}$-functions. In some cases, but not for all, these
functions will result in a function $F(c)$ that satisfies the condition
in Theorem~\ref{main2}. In the examples, we have tended to associate
``lower-order terms'' with the earlier $\Psi$-functions, but what order
is best may require some trial and error.

Whereas we have demonstrated that the main results of the paper are
powerful, regrettably we cannot offer any guarantees that they will
always give results as desired, even when the information matrix can be
written in the form (\ref{infor1}).\vadjust{\goodbreak}

\begin{appendix}\label{append}
\section*{Appendix}

\begin{Proposition}\label{prop0}
Assume that $\{\Psi_0, \Psi_1, \ldots, \Psi_{k-1}\}$ is a Chebyshev
system defined on an interval $[A,B]$. Let $A \leq z_1<z_2<\cdots<z_t
\leq B$, and let $r_1,\ldots,r_t$ be coefficients that satisfy the
following $k$ equations:
\begin{equation}\label{prop0:1}
\sum_{i=1}^{t}r_i\Psi_l(z_{i}) =0,\qquad
l=0,1,\ldots,k-1.
\end{equation}
Then we have:
\begin{longlist}[(b)]
\item[(a)] If $t\leq k$, then $r_i=0, i=1,\ldots,t$.
\item[(b)] If $t=k+1$ and one $r_i$ is not zero, then all are nonzero;
moreover all~$r_i$'s for odd $i$ must then have the same sign, which is
opposite to that of the $r_i$'s for even $i$.
\end{longlist}
\end{Proposition}

\begin{pf}
For part (a), if $t<k$, we can expand $z_1,\ldots,z_t$ to a set of
$k$ distinct points, taking $r_i=0$ for the added points. Thus without
loss of generality, take $t=k$. Consider the matrix
\begin{equation}\label{mat62}
\Psi(z_1,z_2,\ldots,z_k)= \pmatrix{
\Psi_0(z_1) & \Psi_0(z_2) & \cdots& \Psi_0(z_{k}) \cr
\Psi_1(z_1) & \Psi_1(z_2) & \cdots& \Psi_1(z_{k}) \cr
\vdots& \vdots& \ddots& \vdots\cr
\Psi_{k-1}(z_1) & \Psi_{k-1}(z_2) & \cdots& \Psi_{k-1}(z_{k})
}
.
\end{equation}
Then (\ref{prop0:1}) can be written as
\[
\Psi(z_1,z_2,\ldots,z_k) R=0,
\]
where $R=(r_1,\ldots, r_k)^T$. Since $\{\Psi_0, \Psi_1, \ldots,
\Psi
_{k-1}\}$ is a Chebyshev system, $\Psi(z_1,z_2,\ldots,z_k)$ is
nonsingular, so that $R=0$.

For part (b), if one $r_i$ is 0, then it follows from part (a) that
all $r_i$'s are~0. Therefore, if at least one $r_i$ is nonzero, then
all of them must be nonzero. With the notation from the previous
paragraph, we can write (\ref{prop0:1}) as
\[
\Psi(z_1,z_2,\ldots,z_k) R=-r_{k+1} \psi(z_{k+1}),
\]
where $\psi(z_{k+1})=(\Psi_0(z_{k+1}), \Psi_1(z_{k+1}), \ldots,
\Psi
_{k-1}(z_{k+1}))^T$. It follows that
\begin{equation}\label{addthis}
r_i=-r_{k+1}\frac{|\Psi(z_1,\ldots,z_{i-1},z_{k+1},z_{i+1},\ldots
,z_k)|}{|\Psi(z_1,z_2,\ldots,z_k)|},\qquad i=1, \ldots, k.
\end{equation}
By the Chebyshev system assumption, the denominator $|\Psi
(z_1,z_2,\ldots,z_k)|$ in~(\ref{addthis}) is positive, while the
numerator $|\Psi(z_1,\ldots,z_{i-1},z_{k+1},z_{i+1},\ldots,z_k)|$ is
positive for $i = k, k-2, \ldots$ and negative otherwise. The result in
(b) follows.
\end{pf}

\begin{Proposition}\label{prop1}
Let $\{\Psi_0=1, \Psi_1, \ldots, \Psi_{k-1}\}$ be a Chebyshev
system on
an interval $[A,B]$, and suppose that $k=2n-1$. Consider $n$
pairs
$(c_i,\omega_i)$, $i=1,\ldots,n$, and $n$ pairs $(\tilde{c}_i,\tilde
{\omega}_i)$, $i=1,\ldots,n$, with $\omega_i>0$, $\tilde{\omega}_i>0$
and $A\leq c_1<\tilde{c}_1<\cdots<c_n<\tilde{c}_n= B$. Suppose further
that the following\vadjust{\goodbreak} $k$ equations hold:
\begin{equation}\label{prop1:1}
\sum_{i}\omega_i\Psi_l(c_{i})=\sum_{i}\tilde{\omega}_i\Psi
_l(\tilde
{c}_{i}),\qquad l=0,1,\ldots,k-1.
\end{equation}
Then, for any function $\Psi_k$ on $[A,B]$, we can conclude that
\begin{equation}\label{prop1:2}
\sum_{i}\omega_i\Psi_{k}(c_{i})<\sum_{i}\tilde{\omega}_i\Psi
_{k}(\tilde{c}_{i})
\end{equation}
if $\{\Psi_0=1,\Psi_1,\ldots,\Psi_{k-1}, \Psi_k\}$ is also a
Chebyshev system.
\end{Proposition}

\begin{pf}
With
\begin{eqnarray*}
R=(\omega_1, -\tilde{\omega}_1,\omega_2,-\tilde{\omega}_2,\ldots
,\omega_n)^T,
\end{eqnarray*}
the $k$ equations in (\ref{prop1:1}) can be written as
\begin{eqnarray}\label{prop1:3}
\Psi(c_1,\tilde{c}_1,\ldots,c_n) R=\tilde{\omega}_n \psi(\tilde{c}_n),
\end{eqnarray}
where $\Psi$ and $\psi$ are as defined in the proof of Proposition~\ref{prop0}.
Further, (\ref{prop1:2}) is equivalent to
\begin{equation}\label{prop1:4}
(\Psi_k(c_1), \Psi_k(\tilde{c}_1), \ldots, \Psi_k(c_n)) R < \tilde
{\omega}_n \Psi_k(\tilde{c}_n).
\end{equation}
Using (\ref{prop1:3}) to solve for $R$, and using that $\tilde{\omega
}_n > 0$, we see that (\ref{prop1:4}) is equivalent to
\begin{equation}\label{prop1:5}
\qquad(\Psi_k(c_1), \Psi_k(\tilde{c}_1),\ldots, \Psi_k(c_n))\Psi^{-1}(c_1,
\tilde{c}_1,\ldots,c_n) \psi(\tilde{c}_n) - \Psi_k(\tilde{c}_n) < 0.
\end{equation}
From an elementary matrix result [see, e.g., Theorem 13.3.8 of
\citet{Har97}], the left-hand side of (\ref{prop1:5}) can be written as
\begin{equation}\label{prop1:6}
- \frac{|\Psi^*(c_1,\tilde{c}_1,\ldots,c_n,\tilde{c}_n)|}{|\Psi
(c_1,\tilde{c}_1,\ldots,c_n)|},
\end{equation}
where
\begin{eqnarray}\qquad
&&\Psi^*(c_1,\tilde{c}_1,\ldots,c_n,\tilde{c}_n)\nonumber\\[-8pt]\\[-8pt]
&&\qquad= \pmatrix{
\Psi_0(c_1) & \Psi_0(\tilde{c}_1) & \cdots& \Psi_0(c_n) & \Psi
_0(\tilde{c}_n) \cr
\Psi_1(c_1) & \Psi_1(\tilde{c}_1) & \cdots& \Psi_1(c_n) & \Psi
_1(\tilde{c}_n) \cr
\vdots& \vdots& \ddots& \vdots& \vdots\cr
\Psi_{k-1}(c_1) & \Psi_{k-1}(\tilde{c}_1) & \cdots& \Psi_{k-1}(c_n)
&\Psi_{k-1}(\tilde{c}_n)\cr
\Psi_{k}(c_1) & \Psi_{k}(\tilde{c}_1) & \cdots& \Psi_{k}(c_n)
&\Psi
_{k}(\tilde{c}_n)
}.\nonumber
\end{eqnarray}
Since both $\{\Psi_0,\Psi_1,\ldots,\Psi_{k-1}\}$ and $\{\Psi
_0,\Psi
_1,\ldots,\Psi_{k-1}, \Psi_k\}$ are Chebyshev systems and
$c_1<\tilde
{c}_1<\cdots<c_n<\tilde{c}_n$, it follows that (\ref{prop1:6}) is
negative, which is what had to be shown.
\end{pf}

A similar argument as for Proposition~\ref{prop1} can be used for the
next result.

\begin{Proposition}\label{prop2}
Let $\{\Psi_0=1,\Psi_1,\ldots,\Psi_{k-1}\}$ be a Chebyshev system
on an
interval $[A,B]$ and suppose that $k=2n$. Consider $n$ pairs
$(c_i,\omega_i)$, $i=1,\ldots,n$, and $n+1$ pairs\vadjust{\goodbreak}
$(\tilde{c}_i,\tilde{\omega}_i)$, $i=0,1,\ldots,n$, with $\omega_i>0$,
$\tilde{\omega}_i>0$ and $A=\tilde{c}_0< c_1<\tilde{c}_1< \cdots
<c_n<\tilde{c}_n= B$. Suppose further that the following $k$
equations hold:
\begin{equation}\label{prop2:1}
\sum_{i}\omega_i\Psi_l(c_{i})=\sum_{i}\tilde{\omega}_i\Psi
_l(\tilde{c}_{i}),\qquad
l=0,1,\ldots,k-1.
\end{equation}
Then, for any function $\Psi_k$ on $[A,B]$, we can conclude that
\begin{equation}\label{prop2:2}
\sum_{i}\omega_i\Psi_{k}(c_{i})<\sum_{i}\tilde{\omega}_i\Psi
_{k}(\tilde{c}_{i})
\end{equation}
if $\{\Psi_0=1,\Psi_1,\ldots,\Psi_{k-1}, \Psi_k\}$ is also a
Chebyshev system.
\end{Proposition}

\begin{Proposition}\label{prop3}
Consider functions $\Psi_0=1, \Psi_1, \ldots, \Psi_k$ on an interval
$[A,B]$. Compute the corresponding functions $f_{l,l}$ as in (\ref
{def:df}), but with $C_{22}(c)$ replaced by $\Psi_k$, and suppose that
$f_{l,l}>0$, $l=1,\ldots,k-1$. Then
$\{1, \Psi_1, \ldots, \Psi_{k}\}$ is a Chebyshev system if $f_{k,k}>0$,
while $\{1, \Psi_1, \ldots, -\Psi_{k}\}$ is a Chebyshev system if
$f_{k,k}<0$.
\end{Proposition}

\begin{pf}
The conclusion for the case $f_{k,k}<0$ follows immediately from that
for $f_{k,k}>0$, so that we will only focus on the latter. We need to
show that
\begin{equation}\label{prop3:1}
\left|
\matrix{
1 & 1 & \cdots& 1 \cr
\Psi_1(z_0) & \Psi_1(z_1) & \cdots& \Psi_1(z_k) \cr
\vdots& \vdots& \ddots& \vdots\cr
\Psi_{k}(z_0) & \Psi_{k}(z_1) & \cdots& \Psi_{k}(z_k)
}
\right|>0
\end{equation}
for any given $A \leq z_0<z_1<\cdots<z_k\leq B$. Consider (\ref
{prop3:1}) as a function of~$z_k$. The determinant is 0 if $z_k =
z_{k-1}$, so that it suffices to show that the derivative of (\ref
{prop3:1}) with respect to $z_k$ is positive on $(z_{k-1}, B)$, that is,
\begin{equation}\label{prop3:2}
\left|
\matrix{
1 & 1 & \cdots& 1& 0 \cr
\Psi_1(z_0) & \Psi_1(z_1) & \cdots& \Psi_1(z_{k-1}) & f_{1,1}(z_k)
\cr
\vdots& \vdots& \ddots& \vdots& \vdots\cr
\Psi_{k}(z_0) & \Psi_{k}(z_1) & \cdots&\Psi_k(z_{k-1}) & f_{k,1}(z_k)
}
\right|>0
\end{equation}
for any $z_k \in(z_{k-1},B)$. Now consider (\ref{prop3:2}) as a
function of $z_{k-1}$, and use a~similar argument. It suffices to show
that for $z_{k-1}\in(z_{k-2}, z_{k})$,
\begin{equation}
\left|
\matrix{
1 & 1 & \cdots& 0& 0 \cr
\Psi_1(z_0) & \Psi_1(z_1) & \cdots& f_{1,1}(z_{k-1}) & f_{1,1}(z_k)
\cr
\vdots& \vdots& \ddots& \vdots& \vdots\cr
\Psi_{k}(z_0) & \Psi_{k}(z_1) & \cdots&f_{k,1}(z_{k-1}) & f_{k,1}(z_k)
}
\right|>0.
\end{equation}
Continuing like this, it suffices to show that
\begin{equation}\label{prop3:3}
\left|
\matrix{
f_{1,1}(z_1) & f_{1,1}(z_2) & \cdots& f_{1,1}(z_k) \cr
\vdots& \vdots& \ddots& \vdots\cr
f_{k,1}(z_1) & f_{k,1}(z_2) &\cdots& f_{k,1}(z_k)
}
\right|>0\vadjust{\goodbreak}
\end{equation}
for any $A\leq z_1<z_2<\cdots<z_k\leq B$. Since $f_{1,1}(c)>0$ for
$c\in
[A,B]$, (\ref{prop3:3}) is equivalent to
\begin{equation}\label{prop3:4}
\left|
\matrix{
1 & 1 & \cdots& 1 \vspace*{2pt}\cr
\dfrac{f_{2,1}(z_1)}{f_{1,1}(z_1)} & \dfrac{f_{2,1}(z_2)}{f_{1,1}(z_2)}
& \cdots& \dfrac{f_{2,1}(z_k)}{f_{1,1}(z_k)} \cr
\vdots& \vdots& \ddots& \vdots\vspace*{2pt}\cr
\dfrac{f_{k,1}(z_1)}{f_{1,1}(z_1)} & \dfrac{f_{k,1}(z_2)}{f_{1,1}(z_2)}
&\cdots& \dfrac{f_{k,1}(z_k)}{f_{1,1}(z_k)}
}
\right|>0.
\end{equation}
Recall that the entries in the last $k-1$ rows of this matrix are by
definition simply values of $f_{l,2}$, $l=2,\ldots,k$. Hence, applying
the same arguments used for (\ref{prop3:1}) to (\ref{prop3:4}) and
using that $f_{2,2}(c)>0$ for $c\in[A,B]$, it is sufficient to show that
\begin{equation}\label{prop3:5}
\left|
\matrix{
1 & 1 & \cdots& 1 \vspace*{2pt}\cr
\dfrac{f_{3,2}(z_2)}{f_{2,2}(z_2)} & \dfrac{f_{3,2}(z_3)}{f_{2,2}(z_3)}
& \cdots& \dfrac{f_{3,2}(z_k)}{f_{2,2}(z_k)} \cr
\vdots& \vdots& \ddots& \vdots\vspace*{2pt}\cr
\dfrac{f_{k,2}(z_2)}{f_{2,2}(z_2)} & \dfrac{f_{k,2}(z_3)}{f_{2,2}(z_3)}
&\cdots& \dfrac{f_{k,2}(z_k)}{f_{2,2}(z_k)}
}
\right|>0.
\end{equation}
Continuing like this, the ultimate sufficient condition is that
$f_{k,k}(c)>0$ for $c\in[A,B]$, which is precisely our assumption. Thus
the conclusion follows.
\end{pf}
\end{appendix}



\printaddresses

\end{document}